\documentclass[psamsfonts]{amsart}

\usepackage{amssymb,amsfonts}
\usepackage[all,arc]{xy}
\usepackage{enumerate}
\usepackage{mathrsfs}
\usepackage{mathtools}

\usepackage{cite}
\newtheorem{theorem}{Theorem}[section]
\newtheorem{corollary}[theorem]{Corollary}
\newtheorem{proposition}[theorem]{Proposition}

\theoremstyle{definition}
\newtheorem{definition}[theorem]{Definition}

\newtheorem{example}[theorem]{Example}

\newtheorem{remark}[theorem]{Remark}

\theoremstyle{remark}

\DeclareMathOperator{\ZC}{ZC}
\DeclareMathOperator{\Ann}{Ann}
\begin{document}
\title{Eversible and reversible semigroups and semirings}
	
\author{\bfseries Peyman Nasehpour}
	
\address{Department of Engineering Science\\ Golpayegan University of Technology\\ Golpayegan\\
		Iran}
\email{nasehpour@gut.ac.ir, nasehpour@gmail.com}
	
	

\begin{abstract}
The main purpose of this paper is to investigate the zero-divisors of semirings and in particular, to discuss eversible and reversible semirings. Since some of our results can be given in a more general context, we introduce a new ring-like algebraic structure, called prenearsemiring and investigate its left and right zero-divisors as well as its nilpotents. \end{abstract}
	
\maketitle
	
\section{Introduction}

In this paper, we continue our investigations on nilpotents and zero-divisors of some algebraic structures. In order to do that we generalize some concepts such as reversible \cite{Cohn1999} and eversible \cite{GhashghaeiKosanNamdariYildirim2019} rings in the context of semigroups with zero and introduce a new ring-like algebraic structure which we call it prenearsemiring (for some useful resources on ring-like algebraic structures see \cite{Glazek2002,Golan1999(a),Golan1999(b),Golan2003,GondranMinoux2008,HebischWeinert1998,KuichSalomaa1986,Pilz1977} and for the author's papers on zero-divisors check  \cite{EpsteinNasehpour2013, Nasehpour2010,Nasehpour2011,Nasehpour2008,Nasehpour2016,Nasehpour2018,NasehpourPayrovi2010,NasehpourPhD,NasehpourGeorgia}). Here is a brief description of what we do in the current paper:

K\"{o}the's conjecture in ring theory states that if a ring has a non-zero nil right ideal, then it has a non-zero nil ideal. P.M. Cohn calls a ring $R$ reversible if $ab=0$ implies $ba=0$, for all $a,b \in R$ and then proves that in any reversible ring, every nil right ideal is a two-sided nil ideal, and the set of all nilpotent elements is a nil ideal \cite{Cohn1999}. Inspired by this, in Definition \ref{reversibledef}, we define a semigroup with zero $(S,\cdot,0)$ to be reversible if $st=0$ implies $ts=0$, for all $s,t\in S$. Note that $(S,\cdot,0)$ is a semigroup with zero if $(S,\cdot)$ is a semigroup, the set $S$ has at least two elements, and there is a special element $0\in S$, called the zero of $S$, such that $s \cdot 0=0 \cdot s=0$, for all $s\in S$ \cite[p. 2]{Howie1995}. Also, note that if $S$ is semigroup and $I$ a nonempty subset of $S$, then $I$ is a left ideal of $S$ if $SI \subseteq I$, i.e. $sa\in I$, for all $s\in S$ and $a\in I$, and $I$ is a right ideal of $S$ if $IS \subseteq I$. If $I$ is both a left and a right ideal of $S$, then $I$ is called two-sided ideal of $S$. A left, right, or two-sided ideal $I$ of a semigroup $S$ is nil if each of its elements is nilpotent, i.e. for all $s\in I$, we have $s^n=0$, for some natural number $n$ \cite{Howie1976}. In Theorem \ref{CohnreversibleSemigroup}, we generalize Cohn's theorem for reversible rings and prove that in any reversible semigroup with zero, every nil right ideal is a two-sided nil ideal.

Let us recall that a ring $R$ is said to be eversible if every left zero-divisor on $R$ is also a right zero-divisor on $R$ and conversely (i.e. all one-sided zero-divisors are two-sided) \cite{GhashghaeiKosanNamdariYildirim2019}. Inspired by this, in Definition \ref{eversibledef}, we define a semigroup with zero $S$ to be eversible if every left zero-divisor on $S$ is also a right zero-divisor on $S$ and conversely. Here we need to explain that if $S$ is a semigroup with zero, an element $s\in S$ is a left (right) zero-divisor of $S$ if there exists a nonzero $x$ such that $sx = 0$ ($xs = 0$). The set of left (right) zero-divisors of $S$ is denoted by $Z_l(S)$ ($Z_r(S)$). The set of zero-divisors $Z(S)$ of $S$ is $\{s\in S: \exists{~} t \in S-\{0\} (st =0 \vee ts=0) \} = Z_l(S) \cup Z_r(S)$ \cite{Howie1976}. In Example \ref{eversibleEx}, we construct a suitable example to show that this assumption that the set of right zero-divisors of a semigroup with zero is the subset of its set of left zero-divisors does not automatically imply $S$ to be eversible.

In Section \ref{sec:prenearsemiring}, we introduce prenearsemirings. In fact, in Definition \ref{prenearsemiring}, we define an algebra $(S,+,\cdot,0,1)$ of type $(2, 2, 0, 0)$ to be a right prenearsemiring (for short right PN-semiring) if the following conditions hold: 

\begin{enumerate}
	\item $(S,+,0)$ is a unital magma,
	\item $(S,\cdot,1)$ is a monoid,
	\item the right distributive law satisfies, i.e. $(u+v)w=uw+vw$, for all $u,v,w \in S$,
	\item The element 0 is an absorbing element of the monoid $S$, i.e. $s0 = 0s = 0$, for all $s\in S$.
\end{enumerate}

We define left PN-semirings similarly. Note that a set $M$ with a map $M \times M \rightarrow M$, denoted by $(x,y) \mapsto xy$, is a magma \cite[Definition 1.1]{Serre1965}.

We say a PN-semiring is distributive if it is both a left and a right PN-semiring. In Definition \ref{eversibledef}, we define a PN-semiring $S$ to be eversible if every left zero-divisor on $S$ is also a right zero-divisor on $S$ and conversely and in Theorem \ref{eversiblethm1}, we prove that if $S$ is a distributive PN-semiring, then the following statements are equivalent:

\begin{enumerate}
	\item The PN-semiring $S$ is eversible.
	
	\item For any $a, b \in S$, if $ab = 0$ then one of the following occurs:
	
	\begin{enumerate}
		\item There exists $c \in S-\{0\}$ such that $bc = 0 = ca$.
		\item There exists $c \in S-\{0\}$ such that $bc \neq 0$, $ca \neq 0$, and $bca =0$.
	\end{enumerate}
\end{enumerate}

In Section \ref{sec:reversiblesemirings}, we focus on eversible and reversible semirings. Since different authors define semirings differently \cite{Glazek2002}, it is very important to clarify, from the beginning, what we mean by a semiring in this paper. By a semiring, we understand an algebraic structure $(S,+,\cdot,0,1)$ with the following properties:

\begin{enumerate}
	\item $(S,+,0)$ is a commutative monoid,
	\item $(S,\cdot,1)$ is a monoid with $1\neq 0$,
	\item $a(b+c) = ab+ac$ and $(b+c)a = ba+ca$ for all $a,b,c\in S$,
	\item $a\cdot 0 = 0\cdot a = 0$ for all $a\in S$. 
\end{enumerate}

Let us recall that a semiring $S$ is Armendariz if $f=\sum^m_{i=0} a_i X^i$ and $g=\sum^n_{j=0} b_j X^j$ are elements of the polynomial semiring $S[X]$ with $fg=0$, then $a_i b_j = 0$, for all $1 \leq i \leq m$ and $1 \leq j \leq n$ \cite{GuptaKumar2011}. Note that the definition of Armendariz semirings is based on the definition of Armendariz rings defined in Definition 1.1 in the paper \cite{RegeChhawchharia1997}. Also, note that a semiring $S$ is zerosumfree if $s+t =0$ implies $s=t=0$, for all $s,t \in S$ \cite[p. 4]{Golan1999(b)}. In Section \ref{sec:reversiblesemirings}, for a semiring $S$, we prove the following:

\begin{enumerate}
	\item If $S$ is zerosumfree, then $S$ is Armendariz (Proposition \ref{ZSFisArm}).
	
	\item If $S$ is Armendariz, then $S$ is reversible if and only if the polynomial semiring $S[X]$ is reversible (Theorem \ref{reversiblethm1}).
	
	\item If $S$ is zerosumfree, then $S$ is reversible if and only if the power series semiring $S[[X]]$ is reversible (Theorem \ref{reversiblethm2}).
	
	\item $S[X]$ is reversible if and only if the Laurent polynomial semiring $S[X; X^{-1}]$ is reversible (Theorem \ref{ReversibleLaurent}).
\end{enumerate}

Clearly, one of the results of these statements is that if $S$ is a zerosumfree semiring, then the following statements are equivalent:

\begin{enumerate}
	\item $S$ is reversible,
	\item $S[X]$ is reversible,
	\item $S[[X]]$ is reversible,
	\item $S[X; X^{-1}]$ is reversible (see Corollary \ref{ReversibleLaurent}).
\end{enumerate}

In Proposition \ref{expectationreversible1} and Proposition \ref{expectationreversible2}, we construct some other families of reversible (proper) semirings and in Theorem \ref{expectationeversible1}, we show that there are eversible semirings that are not reversible. We recall that a semiring is proper if it is not a ring.

Let $(M,+,0)$ be a commutative additive monoid. The monoid $M$ is said to be a left $S$-semimodule if $S$ is a semiring and there is a function, called scalar product, $\lambda: S \times M \longrightarrow M$, defined by $\lambda (s,m)= sm$ such that the following conditions are satisfied:

\begin{enumerate}
	\item $s(m+n) = sm+sn$ for all $s\in S$ and $m,n \in M$;
	\item $(s+t)m = sm+tm$ and $(st)m = s(tm)$ for all $s,t\in S$ and $m\in M$;
	\item $s\cdot 0=0$ for all $s\in S$ and $0 \cdot m=0$ and $1 \cdot m=m$ for all $m\in M$.
\end{enumerate}

Right semimodules over semirings are defined similarly. If $S$ and $T$ are semirings and $M$ is a left $S$-semimodule and right $T$-semimodule, then we say that $M$ is an $(S,T)$-bisemimodule. For more on semimodules, check Section 14 of Golan's book \cite{Golan1999(b)}.

Let us recall that if $S$ and $T$ are two semirings, $M$ is an $(S,T)$-bisemimodule, and we set $\mathcal{M} =\begin{pmatrix}
S & M \\
0 & T
\end{pmatrix}$ to be the set of all matrices of the form $\begin{pmatrix}
s & m \\
0 & t
\end{pmatrix}$, where $s\in S$, $t\in T$, and $m\in M$, then $\mathcal{M}$ with componentwise addition and the following multiplication

$$\begin{pmatrix}
s_1 & m_1 \\
0 & t_1
\end{pmatrix} \cdot \begin{pmatrix}
s_2 & m_2 \\
0 & t_2
\end{pmatrix} = \begin{pmatrix}
s_1 s_2 & s_1 m_2 + m_1 t_2 \\
0 & t_1 t_2
\end{pmatrix},$$ is a semiring \cite[Example 7.3]{Golan2003}. Such semirings have applications in the automatic parallelization of linear computer codes \cite{BuckerBuschelmanHovland2002}. In Theorem \ref{GExpectationthm3}, we show that if $S$ and $T$ are semirings, $M$ is an $(S,T)$-bisemimodule such that $Z_S(M)=\{0\}$ and $Z_T(M)=\{0\}$, and the semiring $ \mathcal{M}= \begin{pmatrix}
	S & M \\
	0 & T
	\end{pmatrix}$ is eversible, then the semirings $S$ and $T$ are eversible. Note that by $Z_S(M)$ ($Z_T(M)$), we mean the set of zero-divisors of the left $S$-semimodule (right $T$-semimodule) $M$: \[Z_S(M)=\{s\in S : \exists ~ m\in M-\{0\} ~ (sm=0)\} (Z_T(M)=\{t\in T : \exists ~ m\in M-\{0\} ~ (mt=0)\}).\]
	
\section{Eversible and Reversible Semigroups with Zero}\label{sec:reversiblesemigroups}

P.M. Cohn defines a ring $R$ to be reversible if $ab=0$ implies $ba=0$, for all $a,b \in R$ \cite{Cohn1999}. Similarly, reversible semirings have been defined in \cite{GuptaKumar2011}. Based on these, we give the following definition:

\begin{definition}
	
	\label{reversibledef}
	
	A semigroup with zero $(S,\cdot,0)$ is reversible if $st=0$ implies $ts=0$, for all $s,t\in S$.
\end{definition}

\begin{remark}
	Anderson and Camillo \cite{AndersonCamillo1999} used the term rings satisfying $\ZC_2$ for what is called reversible rings by Cohn \cite{Cohn1999}. Krempa and Niewieczerzal in \cite{KrempaNiewieczerzal1977} applied the term $C_0$ for such rings.
\end{remark}

The proof of the following statement is straightforward.

\begin{proposition}
	
	\label{reversiblepro1}
	
	Let $\Lambda$ be an index set and $S_i$ be a reversible semigroup with zero, for each $i\in \Lambda$. Then, $\bigoplus_{i\in \Lambda} S_i$ and $\prod_{i\in \Lambda} S_i$ are reversible semigroups with zero.
\end{proposition}

\begin{definition}
	
	Let $S$ be a semigroup with zero.
	
	\begin{enumerate}
		
		\item We say $S$ is entire if $ab=0$ implies either $a=0$ or $b=0$, for all $a,b\in S$, i.e. $S$ has no proper zero-divisors.
		
		\item We define $S$ to be prime if for any two elements $a,b \in S$, $asb = 0$ for all $s\in S$ implies either $a = 0$ or $b = 0$.
		
		\item  We say $S$ is semiprime if for any $a\in S$, $asa=0$ for all $s\in S$ implies $a=0$.
		
		\item An element $s\in S$ is nilpotent if there is a natural number $n$ such that $s^n =0$. We define $S$ to be nilpotent-free if $s^2=0$ implies $s=0$, for all $s\in S$.
	\end{enumerate}
	
\end{definition}

A multiplicative monoid $(M,\cdot,1)$ is with zero if there is a special element $0\in M$ such that $0\not=1$ and $0m = m0=0$, for all $m\in M$.

\begin{theorem}
	
	\label{reversiblethm3}
	
	The following statements hold:
	
	\begin{enumerate}
		\item A semigroup with zero $(S,\cdot,0)$ is entire if and only if it is prime and reversible;
		\item A monoid with zero $(S,\cdot,0,1)$ is nilpotent-free if and only if it is semiprime and reversible.
	\end{enumerate}
	
	\begin{proof}
		
		(1): Clearly, if $S$ is entire, then $S$ is prime and reversible. Now, let $S$ be prime and reversible and $ab=0$. Obviously, $abs=0$ for each $s\in S$. Since $S$ is reversible, we have that $bsa =0$ for all $s\in S$. At last, since $S$ is prime, we have either $a=0$ or $b=0$.
		
		(2): Let $S$ be nilpotent-free and $ab=0$. So, $(ba)^2 = baba =0$. Since $S$ is nilpotent-free, $ba=0$. This proves the reversibility. For the proof of the semiprimeness of $S$, assume that for an $a\in S$, we have $asa=0$ for all $s\in S$. So, $a1a=0$. Since $S$ is nilpotent-free, $a=0$. Conversely, let $S$ be semiprime and reversible. If $s^2 =0$ for some $s\in S$, then $ssa=0$ for each $a\in S$. Since $S$ is reversible, $sas=0$ for each $a\in S$. Now, since $S$ is semiprime, $s=0$ and the proof is complete.
	\end{proof}
	
\end{theorem}

Let us recall that a ring $R$ is symmetric if $rst = 0$ implies $srt =0$, for all $r$, $s$, and $t$ of $R$. Similarly, we give the following definition:

\begin{definition}
	
	\label{symmetricdef}
	
	A semigroup with zero $(S,\cdot,0)$ is symmetric if $rst = 0$ implies $srt =0$, for all $r,s,t\in S$.
\end{definition}

\begin{theorem}
	
	\label{reversiblethm4}
	
	The following statements hold:
	
	\begin{enumerate}
		\item If a semigroup $(S,\cdot,0)$ is nilpotent-free, then $S$ is symmetric.
		\item If a monoid with zero $(S,\cdot,0,1)$ is symmetric, then $S$ is reversible.
	\end{enumerate}
	
	\begin{proof}
		
		(1): Let $S$ be nilpotent-free. By the proof of Theorem \ref{reversiblethm3}(2), $S$ is reversible. Now, if $ab=0$, then by reversibility, we have that $sba=0$, for each $s\in S$ and this implies that $asb =0$, for each $s\in S$. Now, with the help of this, we prove that $S$ is symmetric. Assume that $rst =0$. So, $rast=0$ and then, $rasbt=0$, for each $a,b \in S$. Using this, we see that $(srt)^2 = s(rtsrt)=0$. This implies that $srt=0$, since $S$ is nilpotent-free. The proof of the statement (2) is straightforward. \end{proof}
\end{theorem}

K\"{o}the's conjecture in ring theory states that if a ring has a non-zero nil right ideal, then it has a non-zero nil ideal. The following is a generalization of Cohn's theorem related to K\"{o}the's conjecture (see Theorem 2.2 in \cite{Cohn1999}):

\begin{theorem}[Cohn's Theorem for Reversible Semigroups] \label{CohnreversibleSemigroup} In any reversible semigroup with zero, every nil right ideal is a two-sided nil ideal. 
	
\begin{proof}
Let $S$ be a reversible semigroup with zero. If $a\in S$ is in such a way that $a^m =0$, then by reversibility, we have $a^{m_1} t_1 a^{m-m_1} =0$, for each $m_1 \leq m$ and $t_1 \in S$. Clearly, this implies that $a^{m_1} t_1 a^{m-m_1} t_2=0$, for each $t_2 \in S$. By induction, it is, then, easy to see that \[a^{m_1} t_1 a^{m_2} t_2 \cdots a^{m_k} t_k =0,\] for arbitrary $t_i \in S$ with arbitrary $m_i$ satisfying $m_1 + m_2 + \dots + m_k = m$. In other words, we have $x_1 x_2 \cdots x_l =0$ as long as $m$ factors of $x_i$ are equal to $a$. By this general argument, we see that if $I$ is a nil right ideal of $S$ and $a\in I$ is in such a way that $a^n=0$ for some positive integer $n$, then $(sa)^n=0$ for each $s\in S$. So, $I$ is two-sided and this completes the proof. \end{proof}
	
\end{theorem}

\begin{definition}
	
	\label{eversibledef}
	
	We define a semigroup with zero $S$ to be eversible if every left zero-divisor on $S$ is also a right zero-divisor on $S$ and conversely, i.e. $Z_l(S) = Z_r(S)$.
	
\end{definition}

\begin{example}
	
	\label{eversibleEx}
	
	Let $S$ be a semigroup with zero. The question may arise if $Z_r(S) \subseteq Z_l(S)$ automatically implies that $Z_l(S) \subseteq Z_r(S)$, i.e. $S$ is eversible. The following example, inspired by the example given in p. 3 of Lam's book on ring theory \cite{Lam2001}, shows that this is not the case. Consider the semigroup $\mathcal{S} = \begin{pmatrix}
	\mathbb N_0 & \mathbb Z_2 \\
	0 & \mathbb Z_2
	\end{pmatrix}$ with formal matrix multiplication. First, we show that each right zero-divisor is a left zero-divisor on $\mathcal{S}$. Let $D = \begin{pmatrix}
	a & b \\
	0 & c
	\end{pmatrix}$ be a nonzero right zero-divisor on $\mathcal{S}$. So, there is a nonzero element $W= \begin{pmatrix}
	x & y \\
	0 & z
	\end{pmatrix}$ in $\mathcal{S}$ such that $ WD = \begin{pmatrix}
	xa & xb+yc \\
	0 & zc
	\end{pmatrix} = \begin{pmatrix}
	0 & 0 \\
	0 & 0
	\end{pmatrix}.$ Clearly, this implies that $x=0$ or $a=0$, and also, $z=0$ or $c=0$ (and, of course, $xb+yc=0$). So, we get four cases:
	
	\begin{description}
		\item[Case (i)] Let $x=0$ and $z=0$. Since $W$ is nonzero, $y\not=0$. But $yc=xb+yc=0$ in $\mathbb Z_2$. Therefore, $c=0$. Now, it is easy to see that the nonzero $D = \begin{pmatrix}
		a & b \\
		0 & 0
		\end{pmatrix}$ is a left zero-divisor on $\mathcal{S}$.
		
		\item[Case (ii)] Let $x=0$ and $z\not=0$. From $zc=0$, we get that $c=0$. Similar to the Case (i), the nonzero $D$ is a left zero-divisor on $\mathcal{S}$.
		
		\item[Case (iii)] Let $x\not=0$ and $z=0$. From $xa=0$, we get that $a=0$. In this case, an easy calculation shows that the nonzero $D = \begin{pmatrix}
		0 & b \\
		0 & c
		\end{pmatrix}$ is also a left zero-divisor on $\mathcal{S}$.
		
		\item[Case (iv)] Let $x\not=0$ and $z\not=0$. In this case, we need to have $a=0$ and $c=0$. Since $D$ is nonzero, $b\not=0$ and obviously, $D$ is a left zero-divisor on $\mathcal{S}$.
	\end{description}
	
	Now, let $A = \begin{pmatrix}
	2 & 0 \\
	0 & 1
	\end{pmatrix}$ and $B = \begin{pmatrix}
	0 & 1 \\
	0 & 0
	\end{pmatrix}$. Clearly, $AB=0$, which means that $A$ is a left zero-divisor on $\mathcal{S}$. Now, we take an arbitrary element $C= \begin{pmatrix}
	x & y \\
	0 & z
	\end{pmatrix}$ in $\mathcal{S}$. A simple calculation shows that $CA = \begin{pmatrix}
	2x & y \\
	0 & z
	\end{pmatrix}$. So, if $CA = 0$, then $C=0$, which shows that $A$ is not a right zero-divisor on $\mathcal{S}$.
\end{example}

Let us recall that if $S$ is a semigroup with zero, a directed graph $\Gamma(S)$, called zero-divisor graph of $S$, is attributed to $S$ whose vertices is the proper zero-divisors of $S$ and $s \rightarrow t$ is an edge of $\Gamma(S)$ between the vertices $s$ and $t$ if $s\not=t$ and $st=0$ \cite{CannonNeuerburgRedmond2005}. We finalize this section by bringing the following result from \cite{CannonNeuerburgRedmond2005,Redmond2002}, though written in our terminology (eversible semigroups):

\begin{theorem}\cite[Theorem 2.1]{CannonNeuerburgRedmond2005} Let $S$ be a semigroup with zero. The directed graph $\Gamma(S)$ is connected if and only if $S$ is eversible. Moreover, if $\Gamma(S)$ is connected, then the diameter of the graph $\Gamma(S)$ is at most 3.
\end{theorem}

Now we pass to the next section to discuss zero-divisors of those multiplicative semigroups with zero that they also have some kind of additive structure, i.e. they are ring-like algebraic structures. 

\section{Eversible and Reversible Prenearsemirings}\label{sec:prenearsemiring}

In this section, first we introduce a new ring-like algebraic structure:

\begin{definition}
	
	\label{prenearsemiring}
	
We define an algebra $(S,+,\cdot,0,1)$ of type $(2, 2, 0, 0)$ to be a right prenearsemiring (for short PN-semiring) if the following conditions hold: 
	
	\begin{enumerate}
		\item $(S,+,0)$ is a unital magma, i.e. 0 is the neutral element of the magma $(S,+)$,
		\item $(S,\cdot,1)$ is a monoid,
		\item the right distributive law satisfies, i.e. $(u+v)w=uw+vw$, for all $u,v,w \in S$,
		\item the element 0 is an absorbing element of the monoid $S$, i.e. $s0 = 0s = 0$, for all $s\in S$.
	\end{enumerate}

A left PN-semiring is an algebra $(S,+,\cdot,0,1)$ of type $(2, 2, 0, 0)$ satisfying the conditions 1, 2, 4, and the left distributive law, i.e. $u(v+w)=uv+uw$, for all $u,v,w \in S$, instead of the previous condition 3. A PN-semiring is distributive if it is both a left and a right PN-semiring. In this paper, by a PN-semiring, we mean a right PN-semiring if not explicitly stated otherwise.
\end{definition}

In the following, we give a general example for PN-semirings: 

\begin{example}
	
	\label{PNsemiringEx}
	
	The French mathematician Jean-Pierre Serre calls a set $M$ with a map $M \times M \rightarrow M$ denoted by $(x,y) \mapsto xy$, a magma \cite[Definition 1.1]{Serre1965}. A morphism of magmas is a function $f: M \rightarrow N$ mapping from the magma $M$ to the magma $N$, that preserves the binary operation, i.e. $f(xy) = f(x) f(y)$, for all $x,y \in M$. A morphism $f: M \rightarrow M$ is a magma endomorphism of $M$. A ``unital magma endomorphism'' is a magma endomorphism $f:M \rightarrow M$ such that it also preserves the neutral element, i.e. $f(0) = 0$. 
	
	Now, let $(M,+,0)$ be a unital magma. We gather all unital magma endomorphisms of $M$ in the set $E_0(M)$ and define the following binary operations:
	
	\begin{enumerate}
		\item[A:] $(f+g)(m)=f(m)+g(m)$, for all $m\in M$ and $f,g \in E_0(M)$,
		\item[M:] $(f \circ g)(m) = f(g(m))$, for all $m\in M$ and $f,g \in E_0(M)$.
	\end{enumerate}
	
	Then it is easy to verify that $(E_0(M),+,\circ,0,\iota_M)$ is a PN-semiring, where $0: M \rightarrow M$ defined by $0(m)=0$ is the neutral element of the operation + and the identity map $\iota_M$ defined on $M$ by $\iota_M(m)=m$ is the neutral element of the operation $\circ$.
\end{example}

\begin{definition}
	We define a PN-semiring $(S,+,\cdot,0,1)$ to be eversible (reversible) if $(S,\cdot,0)$ as a semigroup with zero is eversible (reversible).
\end{definition}

\begin{theorem}
	
	\label{eversiblethm1}
	
	Let $S$ be a distributive PN-semiring. Then the following statements are equivalent:
	
	\begin{enumerate}
		\item The PN-semiring $S$ is eversible.
		
		\item For any $a, b \in S$, if $ab = 0$ then one of the following occurs:
		
		\begin{enumerate}
			\item There exists $c \in S-\{0\}$ such that $bc = 0 = ca$.
			\item There exists $c \in S-\{0\}$ such that $bc \neq 0$, $ca \neq 0$, and $bca =0$.
		\end{enumerate}
	\end{enumerate}
	
	\begin{proof}
		$(1) \Rightarrow (2)$: Let $ab = 0$ and there is no $c \neq 0$ such that $bc = 0 = ca$. Since $S$ is eversible, there exist $d,e \neq 0$ such that $bd = 0 = ea$ while $be \neq 0$ and $0 \neq da$. Our claim is that $b(d + e) \neq 0$ and
		$(d + e)a \neq 0$. In contrary, if $ b(d + e)=0$ then $bd+be =0$, which implies that $be=ea=0$, a contradiction. Similarly, $(d + e)a \neq 0$. Clearly, $b(d + e)a = 0$.
		
		$(2) \Rightarrow (1)$: Let $x$ be a left zero-divisor. So, there exists $y \neq 0$ such that
		$xy = 0$. If there exists $c \neq 0$ such that $yc = 0 = cx$ then it is clear that $x$ is
		a right zero-divisor. Otherwise, there exists $c \neq 0$ such that $yc \neq 0$, $cx \neq 0$, and
		$ycx = 0$. Since $yc \neq 0$, we conclude that $x$ is a right zero-divisor. A dual proof shows that any right zero-divisor is also a left zero-divisor. So, $S$ is eversible and this completes the proof. \end{proof}
	
\end{theorem}

\begin{example}
	Let us recall that $(S,+,\cdot,0,1)$ is a left (right) seminearring with an identity element 1 if $(S,+,0)$ and $(S,\cdot,1)$ are monoids and left (right) distributive law holds and 0 is an absorbing element of the multiplicative monoid $(S,\cdot)$ \cite[Definition 1]{vanHoornvanRootselaar1967}. Note that almost the same algebraic structure with the German name ``Fasthalbring'' (translated into English as nearsemiring) was introduced by van Rootselaar in \cite{vanRootselaar1962}. Now, it is clear that examples of distributive (i.e. both left and right) PN-semirings include distributive seminearrings with identity, distributive nearrings with identity \cite{Pilz1977}, semirings \cite{Golan1999(b)}, and of course rings \cite{Lam2001}.
\end{example}

Inspired by the definition of the ideals in the semiring theory \cite{Bourne1951}, we give the following definition:

\begin{definition} 
	
	\label{PNsemiringideal}
	
	Let $S$ be a right PN-semiring and $I$ a nonempty subset of $S$.
	
	\begin{enumerate}
		\item We say $I$ is a left ideal of $S$ if it is an additive submagma of $S$ and absorbing multiplication on the left, i.e. $a+b\in I$ and $sa\in I$, for all $a,b\in I$ and $s\in S$. We define the right ideals similarly. We say $I$ is a two-sided ideal of $S$ if it is both a left and a right ideal of $S$.
		
		\item We call a left, right, or two-sided ideal of the PN-semiring $S$ nil if each of its elements is nilpotent.
	\end{enumerate}

The above concepts for left PN-semirings can be defined similarly.
 
\end{definition}

\begin{theorem}[Cohn's Theorem for Reversible PN-Semirings] \label{Cohnreversible1} In any reversible PN-semiring, every nil right ideal is a two-sided nil ideal. 
	
	\begin{proof}
		The proof is the same as the proof of Theorem \ref{CohnreversibleSemigroup} and so, omitted.\end{proof}
	
\end{theorem}

Now, we pass to the next section to discuss nilpotent and zero-divisors of semirings.

\section{Eversible and Reversible Semirings}\label{sec:reversiblesemirings}

\begin{theorem}[Cohn's Theorem for Reversible Semirings] \label{Cohnreversible2} In any reversible semiring, every nil right ideal is a two-sided nil ideal, and the set of all nilpotent elements is a nil ideal.
	
\begin{proof}
	
By Theorem \ref{Cohnreversible1}, every nil right ideal of a semiring is a two-sided nil ideal. Now, let $N$ be the set of all nilpotent elements of $S$. If $s$ and $t$ are both in $N$, then there are two positive integers $m$ and $n$ such that $s^m = t^n =0$. It is clear that $(s+t)^{m+n-1} =0$, which means that $s+t \in N$. On the other hand, if $s\in N$, then $(asb)^m=0$, for all $a$ and $b$ of $S$. So, $N$ is a two-sided nil ideal and the proof is complete.
\end{proof}

\end{theorem}

\begin{proposition}
	The class of reversible semirings is closed under subrings and direct products.
	
	\begin{proof}
		Straightforward.\end{proof}
\end{proposition}

\begin{proposition}
	
	\label{ZSFisArm}
	
	Each zerosumfree semiring is Armendariz.
	
	\begin{proof} Let $S$ be zerosumfree and $f,g \in S[X]$ be polynomials such that $f=\sum^m_{i=0} a_i X^i$, $g=\sum^n_{j=0} b_j X^j$ and $fg=0$. Therefore, $c_k=a_0 b_k + a_1 b_{k-1} + \cdots + a_k b_0 =0$, for arbitrary $k$, where $c_k$s are the coefficients of the polynomial $fg$. Since $S$ is zerosumfree, we have $a_i b_j=0$ for all $i,j$. This shows that $S$ is Armendariz. \end{proof}
\end{proposition}

\begin{theorem}
	
	\label{reversiblethm1}
	
	Let $S$ be an Armendariz semiring. Then the semiring $S$ is reversible if and only if the polynomial semiring $S[X]$ is reversible.

\begin{proof} It is straightforward to see that if $S[X]$ is reversible, then so is the semiring $S$. Now, let $S$ be an Armendariz and a reversible semiring and $f$ and $g$ be two polynomials over $S$ such that $fg=0$. Since $S$ is Armendariz, we have $a_i b_j=0$, for all $i,j$. Since $S$ is reversible, $b_j a_i=0$, for all $i,j$. This, obviously, implies that $gf=0$. \end{proof} \end{theorem}

\begin{theorem}
	
	\label{reversiblethm2}
	
	Let $S$ be a zerosumfree semiring. Then the semiring $S$ is reversible if and only if the power series semiring $S[[X]]$ is reversible.
	
	\begin{proof} The proof is similar to the proof of Theorem \ref{reversiblethm2}, and so omitted. \end{proof}
	
\end{theorem}

\begin{corollary}
	
	\label{reversiblecor}
	
	Let $S$ be a zerosumfree semiring. Then the following statements are equivalent:
	
	\begin{enumerate}
		\item $S$ is reversible,
		\item $S[X]$ is reversible,
		\item $S[[X]]$ is reversible.
	\end{enumerate}
\end{corollary}

\begin{example} The semiring $\mathcal{M}_n(S)$ of all $n$-by-$n$ matrices over each arbitrary semiring $S$ is never reversible.
	
	\begin{proof} For $\mathcal{M}_n(S)$, define $A=(a_{ij})$ and $B=(b_{ij})$ in $\mathcal{M}_n(S)$ as follows: \[a_{11} =b_{n1} = b_{n2} = \dots = b_{nn} =1\] and the other arrays of $A$ and $B$ are all zero. Now it is clear that $A B = 0$, while $BA \neq 0$. For example, in the case $n=2$, we have the following:
	\[\begin{pmatrix}
	1 & 0 \\
	0 & 0
	\end{pmatrix} \cdot \begin{pmatrix}
		0 & 0 \\
		1 & 1 \end{pmatrix} = \begin{pmatrix}
		0 & 0 \\
		0 & 0
	\end{pmatrix}, \text{while} \begin{pmatrix}
	0 & 0 \\
	1 & 1
	\end{pmatrix} \cdot \begin{pmatrix}
	1 & 0 \\
	0 & 0 \end{pmatrix} = \begin{pmatrix}
	0 & 0 \\
	1 & 0
	\end{pmatrix}.\] \end{proof}
\end{example}

Let $S$ be a semiring and $_S M_S$ an $(S,S)$-bisemimodule. On the set $S \times M$, we define the two addition and multiplication operations as follows:

\begin{enumerate}
	\item[A:] $(s_1, m_1) + (s_2,m_2) = (s_1 + s_2 , m_1 + m_2)$,
	
	\item[M:] $(s_1, m_1) \cdot (s_2, m_2) = (s_1 s_2, s_1 m_2 + m_1 s_2)$.
\end{enumerate}

The set $S \times M$, equipped with the above operations, denoted by $S \widetilde{\oplus} M$, is a semiring called the expectation semiring of the bisemimodule $M$. For more on expectation semirings in commutative semiring theory, check \cite{NasehpourExpectation}.

\begin{proposition}
	
	\label{expectationreversible1}
	
	Let $S$ be a nilpotent-free semiring. The expectation semiring $S \widetilde{\oplus} S$ is reversible.
	
	\begin{proof} Let $S$ be a nilpotent-free semiring. Consider two arbitrary elements $(a,b)$ and $(c,d)$ in $S \widetilde{\oplus} S$ such that $(a,b)(c,d)=(0,0)$. Clearly, this implies that $ac=0$ and $ad+bc=0$. But $S$ is nilpotent-free and so reversible, by Theorem \ref{reversiblethm3}. Therefore, $ca=0$. On the other hand, $0=cad+cbc=cbc$, which implies $(bc)^2 = 0$. Since $S$ is nilpotent-free, we have $bc=0$. So, $ad = 0$ and finally, $cb = 0 = da$. Hence, $(c,d)(a,b)=(0,0)$ and the proof is complete.\end{proof}
\end{proposition}

If $S$ is a semiring, by a semiring endomorphism on $S$, it is meant a semiring homomorphism from $S$ to itself \cite[p. 159]{Golan1999(b)}. Let $S$ be a commutative semiring, $M$ an $S$-semimodule, and $\sigma$ an endomorphism on $S$. One can give $S \times M$ a semiring structure with pointwise addition and the following multiplication:
\[(s_1, m_1) \cdot (s_2, m_2) = (s_1 s_2, \sigma(s_1) m_2 + m_1 s_2).\]

We call this semiring, denoted by $(S \widetilde{\oplus} M)_\sigma$, the $\sigma$-expectation semiring of the $S$-semimodule $M$.

\begin{proposition}
	
\label{expectationreversible2}
	
Let $S$ be a commutative entire semiring and $\sigma$ an endomorphism on $S$ with $\ker(\sigma)=\{0\}$. Then $(S \widetilde{\oplus} S)_\sigma$ is reversible.

\begin{proof} Let $(a,b)$ and $(c,d)$ be elements of $(S \widetilde{\oplus} S)_\sigma$ such that $(a,b)(c,d) = (0,0)$. Clearly, this implies that $ac=0$ and $\sigma(a)d+bc=0$. Since $S$ is entire, either $a=0$ or $c=0$.

Case 1: If $a=0$, then $bc=0$. Obviously, this implies that either $b=0$ or $c=0$. So, $(c,d)(a,b)=(ca,\sigma(c)b+da) =(0,0)$. 

Case 2: If $c=0$, then $\sigma(a)d=0$. Obviously, this implies either $a=0$ or $d=0$, since $\ker(\sigma)=\{0\}$. Therefore, also in this case, we have $(c,d)(a,b)=(ca,\sigma(c)b+da) =(0,0)$. \end{proof} 
\end{proposition}

Obviously, each reversible semiring is eversible but there are some eversible semirings which are not reversible, as we show in the following:

\begin{theorem}
	
	\label{expectationeversible1}
	
	Let $(E,+,\cdot)$ be a commutative entire semiring. Let $S=E \oplus E$ be a semiring with pointwise addition and multiplication. Define $\sigma$ on $S$ by $\sigma(s,t) = (t,s)$. Then the following statements hold:
	
	\begin{enumerate}
		\item The semiring $S$ is nilpotent-free and $\sigma$ is an endomorphism on $S$.
		\item The $\sigma$-expectation semiring $(S \widetilde{\oplus} S)_\sigma$ is eversible.
		\item \label{notreversible} If $1+1 \neq 0$ in $E$, then $(S \widetilde{\oplus} S)_\sigma$ is not reversible.
	\end{enumerate}
	
	\begin{proof} The proof of the statement (1) is straightforward. (2): In order to show that the semiring $(S \widetilde{\oplus} S)_\sigma$ is eversible, we prove that $((a,b),(c,d))$ is a left (right) zero-divisor if and only if either $a = 0$ or $b = 0$.
	
	 First, suppose that $((a, b), (c, d))$ is a left zero-divisor while $a \neq 0$ and
	$b \neq 0$. Take a nonzero element $((e, f ), (g, h))$ such that \[((a, b), (c, d))((e, f ), (g, h)) = 0.\] Therefore, $((ae, bf ), \sigma((a, b))(g, h) + (c, d)(e, f )) = 0$ and so, $((ae, bf ), (bg + ce, ah
	+ df )) = 0$. This means that $ae = bf = bg+ce = ah+df = 0$. Since $E$ is entire, we have
	$e = f = 0$. So, $bg = ah = 0$ and $g = h = 0$. This implies that $((e, f ), (g, h)) = 0$, a
	contradiction. Conversely, suppose that $((a, b), (c, d))$ is an arbitrary element of $(S \widetilde{\oplus} S)_\sigma$ in which either $a = 0$ or $b = 0$. If $a = 0$ then $((0, b), (c, d))((0, 0), (0, 1)) = 0$.
	Also, if $b = 0$ then $((a, 0), (c, d))((0, 0), (1, 0)) = 0$. Hence, $((a, b), (c, d))$ is a
	left zero-divisor.
	
	Now, suppose that $((a, b), (c, d))$ is a right zero-divisor while $a \neq 0$ and
	$b  \neq 0$. Take a nonzero element $((e, f ), (g, h))$ such that \[((e, f ), (g, h))((a, b), (c, d)) = 0.\] So, $((ea, f b), \sigma((e, f ))(a, b) + (g, h)(c, d)) = 0$ and therefore, $((ea, f b), (f c + ga, ed + hb)) = 0$. Consequently, $ea = f b = f c + ga = ed + hb = 0$. Since $E$ is entire, $e = f = 0$. So, $hb = ga = 0$ and $g = h = 0$. This means that $((e, f ), (g, h)) = 0$, a
	contradiction. Conversely, suppose that $((a, b), (c, d))$ is an arbitrary element of $(S \widetilde{\oplus} S)_\sigma$
	where either $a = 0$ or $b = 0$. If $a = 0$ then  $((0, 0), (1, 0))((0, b), (c, d)) = 0$.
	Also, if $b = 0$ then $((0, 0), (0, 1))((0, b), (c, d)) = 0$. Hence, ((a, b), (c, d)) is a
	right zero-divisor and this finishes the proof of (2).
	
	(3): Observe the following computations:
	\[((0,1),(0,1))\cdot ((1,0),(0,1)) = ((0,0), \sigma((0,1))(0,1) + (1,0)(0,1)) = 0,\] while 
	$((1,0),(0,1))\cdot ((0,1),(0,1)) = ((0,0), \sigma((1,0))(0,1) + (0,1)(0,1)) = ((0,0),(0,1+1)) \neq 0.$\end{proof}
	
\end{theorem}

\begin{remark}
	\label{eversiblebutnotreversible} The statement \ref{notreversible} in Theorem \ref{expectationeversible1} is inspired by Example 1.15 in \cite{KimLee2003}.
\end{remark}

Similar to ring theory, we say an element $s$ of a semiring $S$ is regular if it is not a left or a right zero-divisor. We say an element $c$ of a semiring $S$ is central if $cs = sc$, for all $s\in S$. For the concepts of left {\O}re sets and classical left semiring of fractions, see Chapter 11 in Golan's book \cite{Golan1999(b)}.

\begin{theorem}
	
	\label{ClassicalLeftSemiringofFractions}
	
	Let $S$ be a semiring and $O$ be a left {\O}re set in $S$ consisting of central regular elements. Then the classical left semiring of fractions $O^{-1}S$ is reversible (eversible) if and only if the semiring $S$ is reversible (eversible).
	
	\begin{proof}
	Reversible Case: Let $S$ be a reversible semiring and $u^{-1}s \cdot v^{-1}t =0$, where $u,v\in O$ and $s,t\in S$. Since $u$ and $v$ are central, we have $(uv)^{-1} st = u^{-1}s \cdot v^{-1}t =0$. Since $u$ and $v$ are regular, we have $st =0$. Finally, by reversibility, we have $ts=0$. This implies that $v^{-1}t \cdot u^{-1}s =0$. So, $O^{-1}S$ is reversible. It is straightforward to see that if $O^{-1}S$ is reversible, then so is the semiring $S$.
	
	Eversible Case: Assume that $S$ is an eversible semiring and $u^{-1}s \in O^{-1}S$ is a left zero-divisor in $O^{-1}S$, where $u\in O$ and $s\in S$. By definition, there exists a nonzero $v^{-1}t$ in $O^{-1}S$ such that $u^{-1}s \cdot v^{-1}t = 0$. Since $uv$ is invertible, $st=0$. So, there exists a nonzero $z$ in $S$ such that $zs=0$. This implies that $v^{-1}z \cdot u^{-1}s =0$, which means that $u^{-1}s$ is a right zero-divisor. It is straightforward to see that if $O^{-1}S$ is eversible, then so is the semiring $S$.\end{proof}
\end{theorem}

Let $S$ be a semiring. By $S[X; X^{-1}]$, we mean the Laurent polynomial semiring over $S$.

\begin{theorem}
	
	\label{ReversibleLaurent}
	
	Let $S$ be a semiring. Then the following statements hold:
	
	\begin{enumerate}
		\item $S[X]$ is reversible if and only if $S[X; X^{-1}]$ is reversible.
		
		\item $S[X]$ is eversible if and only if $S[X; X^{-1}]$ is eversible.
	\end{enumerate}

\begin{proof} (1): Suppose that $S[X]$ is reversible. Let $f(X)$ and $g(X)$ be Laurent polynomials over the semiring $S$ such that $f(X) g(X) =0$. Clearly, there is a non-negative integer $n$ such that $f_1(X) = X^n f(X)$ and $g_1(X) = X^n g(X)$ are elements of $S[X]$. It is, now, obvious that $f_1(X)g_1(X) =0$ and so, $g_1(X)f_1(X) =0$. This, obviously, implies that $g(X)f(X) X^{-2n}=0$ and this implies that $g(X)f(X)=0$. Hence, we have already proved that $S[X; X^{-1}]$ is reversible. It is obvious that if $S[X; X^{-1}]$ is reversible, then $S[X]$ is also reversible.
	
(2): The proof of this statement is similar to the proof of the statement (1) and therefore, omitted. \end{proof}
\end{theorem}

\begin{remark}
	Theorem \ref{ReversibleLaurent} is, in fact, a corollary of \ref{ClassicalLeftSemiringofFractions}, since $S[X; X^{-1}] = O^{-1} S[X]$, where $O=\{1,X,X^2,\dots,X^n,\dots\}$, but we gave a direct proof for Theorem  \ref{ReversibleLaurent} in the above.
\end{remark}

\begin{corollary}
		
		\label{reversiblecor2}
		
		Let $S$ be a zerosumfree semiring. Then the following statements are equivalent:
		
		\begin{enumerate}
			\item $S$ is reversible,
			\item $S[X]$ is reversible,
			\item $S[X; X^{-1}]$ is reversible,
			\item $S[[X]]$ is reversible.
		\end{enumerate}
	
	\begin{proof} Combine Corollary \ref{reversiblecor} and Theorem \ref{ReversibleLaurent}.\end{proof}
	
	\end{corollary}

Similar to ring theory \cite[Lemma 10.29]{Lam2001}, we define left and right annihilators of an element of a semiring as follows:

\begin{definition}
	
	\label{leftrightannihilatorsdef}
	
Let $S$ be a semiring and $s\in S$. 

\begin{enumerate}
	\item[(L)] We define the left annihilators of the element $s$ to be the following set: \[\Ann_l(s) = \{x\in S : xs=0\}.\] 
	
	\item[(R)] We define the right annihilators of the element $s$ to be the following set: \[\Ann_r(s) = \{x\in S : sx=0\}.\]
\end{enumerate}

\end{definition}

\begin{theorem}
	
	\label{GExpectationthm1}
	
	Let $S$ and $T$ be semirings and $M$ be an $(S,T)$-bisemimodule and $m$ be a nonzero element of $M$. Then the following statements are equivalent:
	
	\begin{enumerate}
		\item The element $X = \begin{pmatrix}
			s & m \\
			0 & t
		\end{pmatrix}$ of the semiring $ \begin{pmatrix}
		S & M \\
		0 & T
		\end{pmatrix}$ is a left zero-divisor,
		
		\item One of the following conditions are satisfied:
		
		\begin{enumerate}
			\item $s\in Z_l(S)$.
			\item $s\in Z_l(S)-\{0\}$, $t\in Z_l(T)$, and there exists a nonzero $t^{\prime\prime} \in \Ann_r(t)$ such that $m t^{\prime\prime} =0$.
			\item $s\in Z_l(S)-\{0\}$, $t\in Z_l(T)-\{0\}$, and there exists a nonzero $m^{\prime\prime} \in M$ such that $s m^{\prime\prime} =0$.
		\end{enumerate}
	\end{enumerate}

\begin{proof}
	
$(1) \Rightarrow (2)$: Let $X = \begin{pmatrix}
s & m \\
0 & t
\end{pmatrix}$ be a left zero-divisor. So, there exists a nonzero $X^{\prime} = \begin{pmatrix}
s^{\prime} & m^{\prime} \\
0 & t^{\prime}
\end{pmatrix}$ such that $XX^{\prime} = \begin{pmatrix}
ss^{\prime} & sm^{\prime} + mt^{\prime} \\
0 & tt^{\prime}
\end{pmatrix} = 0.$ Now, if $s\notin Z_l(S)$, then $s^{\prime} =0$. And if $t\notin Z_l(T)$, then $t^{\prime} =0$. Finally, if there is no nonzero $m^{\prime\prime}$ in $M$ such that $sm^{\prime\prime}=0$, we have $m^{\prime}=0$, a contradiction.

$(2) \Rightarrow (1)$: Let $X = \begin{pmatrix}
s & m \\
0 & t
\end{pmatrix}$ be an element of the semiring $ \begin{pmatrix}
S & M \\
0 & T
\end{pmatrix}.$

\begin{itemize}
	\item If (a) holds, then define $X^{\prime} = \begin{pmatrix}
	s^{\prime} & 0 \\
	0 & 0
	\end{pmatrix}$, where $s^{\prime}\in \Ann_r(s)$.
	\item If (b) holds, then define $X^{\prime} = \begin{pmatrix}
	0 & 0 \\
	0 & t^{\prime\prime}
	\end{pmatrix}.$
	
	\item If (c) holds, then define $X^{\prime} = \begin{pmatrix}
	0 & m^{\prime\prime} \\
	0 & 0
	\end{pmatrix}.$
\end{itemize}

Now, it is clear that in each case, we have $XX^{\prime} =0$ which means that $X$ is a left zero-divisor. \end{proof}

\end{theorem}

\begin{theorem}
	
	\label{GExpectationthm2}
	
	Let $S$ and $T$ be semirings and $M$ be an $(S,T)$-bisemimodule and $m$ be a nonzero element of $M$. Then the following statements are equivalent:
	
	\begin{enumerate}
		\item The element $X = \begin{pmatrix}
		s & m \\
		0 & t
		\end{pmatrix}$ of the semiring $ \begin{pmatrix}
		S & M \\
		0 & T
		\end{pmatrix}$ is a right zero-divisor,
		
		\item One of the following conditions are satisfied:
		
		\begin{enumerate}
			\item $s\in Z_r(S)$.
			\item $s\in Z_r(S)-\{0\}$, $t\in Z_r(T)$, and there exists a nonzero $t^{\prime\prime} \in \Ann_l(t)$ such that $m t^{\prime\prime} =0$.
			\item $s\in Z_r(S)-\{0\}$, $t\in Z_r(T)-\{0\}$, and there exists a nonzero $m^{\prime\prime} \in M$ such that $s m^{\prime\prime} =0$.
		\end{enumerate}
	\end{enumerate}

\begin{proof} The proof is similar to the proof of Theorem \ref{GExpectationthm1} and so, omitted.\end{proof}
 
\end{theorem}

\begin{theorem}

\label{GExpectationthm3}
	
Let $S$ and $T$ be semirings and $M$ be an $(S,T)$-bisemimodule such that $Z_S(M)=\{0\}$ and $Z_T(M)=\{0\}$. If the semiring $ \mathcal{M}= \begin{pmatrix}
S & M \\
0 & T
\end{pmatrix}$ is eversible, then the semirings $S$ and $T$ are eversible.

\begin{proof}
	
	(L): Let $s\neq 0$ be a left zero-divisor of the semiring $S$. Clearly, $X = \begin{pmatrix}
	s & 0 \\
	0 & 1
	\end{pmatrix}$ is a left zero-divisor of the semiring $ \mathcal{M} = \begin{pmatrix}
	S & M \\
	0 & T
	\end{pmatrix}.$ Since by assumption, $\mathcal{M}$ is eversible, $X$ is a right zero-divisor of $\mathcal{M}$, which means that there is a nonzero $X^{\prime} = \begin{pmatrix}
	s^{\prime} & m^{\prime} \\
	0 & t^{\prime}
	\end{pmatrix}$ in $\mathcal{M}$ such that $X^{\prime} X =0$. This implies that $m^{\prime}$ and $t^{\prime}$ are zero and $s^{\prime} s=0$ while $s^{\prime}\neq 0$. So, $s$ is a right zero-divisor of $S$.
(R): Now, let $s\neq 0$ be a right zero-divisor of $S$. It is clear that $X$ is a right zero-divisor also. Since $\mathcal{M}$ is eversible, there is a nonzero $X^{\prime} = \begin{pmatrix}
	s^{\prime} & m^{\prime} \\
	0 & t^{\prime}
	\end{pmatrix}$ in $\mathcal{M}$ such that $XX^{\prime} =0$. It is clear that this implies that $ss^{\prime}=0$, $t^{\prime}=0$, and $sm^{\prime}=0$. By assumption, $Z_S(M)=\{0\}$. So, $m^{\prime}=0$. Since $X^{\prime} \neq 0$, we have $s^{\prime}\neq 0$. So, $s$ is a left zero-divisor. Hence, $S$ is eversible.
	
	Similarly, with the help of the assumption $Z_T(M)=\{0\}$, it is possible to prove that $T$ is eversible and the proof is complete. \end{proof}

\end{theorem}
	
\subsection*{Acknowledgments}

The author is supported by the Department of Engineering Science at the Golpayegan University of Technology and his special thanks go to the Department for providing all necessary facilities available to him for successfully conducting this research. The author is also grateful to the President of the Golpayegan University of Technology, Professor Dara Moazzami, for his support and encouragements.

\end{document}